\documentclass[11pt]{amsart}
\usepackage[parfill]{parskip} 
\usepackage{amsmath, amsxtra, amsthm, amssymb,eucal,dsfont}
\usepackage[all]{xy}
\usepackage{color}
\usepackage[colorlinks=true, pdfstartview=FitV, linkcolor=blue, citecolor=blue, urlcolor=blue]{hyperref}

\usepackage{graphicx}

\numberwithin{equation}{section}
\swapnumbers

\newcommand{\nc}{\newcommand}
\nc{\rc}{\renewcommand}

\rc{\b}{\mathbb}
\rc{\c}{\mathcal}
\nc{\on}{\operatorname}
\nc{\tn}{\textnormal}

\nc{\bA}{\b A}
\nc{\bB}{\b B}
\nc{\bC}{\b C}
\nc{\bD}{\b D}
\nc{\bE}{\b E}
\nc{\bF}{\b F}
\nc{\bG}{\b G}
\nc{\bH}{\b H}
\nc{\bI}{\b I}
\nc{\bJ}{\b J}
\nc{\bK}{\b K}
\nc{\bL}{\b L}
\nc{\bM}{\b M}
\nc{\bN}{\b N}
\nc{\bO}{\b O}
\nc{\bP}{\b P}
\nc{\bQ}{\b Q}
\nc{\bR}{\b R}
\nc{\bS}{\b S}
\nc{\bT}{\b T}
\nc{\bU}{\b U}
\nc{\bV}{\b V}
\nc{\bW}{\b W}
\nc{\bX}{\b X}
\nc{\bY}{\b Y}
\nc{\bZ}{\b Z}

\nc{\cA}{\c A}
\nc{\cB}{\c B}
\nc{\cC}{\c C}
\nc{\cD}{\c D}
\nc{\cE}{\c E}
\nc{\cF}{\c F}
\nc{\cG}{\c G}
\nc{\cH}{\c H}
\nc{\cI}{\c I}
\nc{\cJ}{\c J}
\nc{\cK}{\c K}
\nc{\cL}{\c L}
\nc{\cM}{\c M}
\nc{\cN}{\c N}
\nc{\cO}{\c O}
\nc{\cP}{\c P}
\nc{\cQ}{\c Q}
\nc{\cR}{\c R}
\nc{\cS}{\c S}
\nc{\cT}{\c T}
\nc{\cU}{\c U}
\nc{\cV}{\c V}
\nc{\cW}{\c W}
\nc{\cX}{\c X}
\nc{\cY}{\c Y}
\nc{\cZ}{\c Z}

\nc{\fA}{{\mathfrak A}}
\nc{\fB}{{\mathfrak B}}
\nc{\fC}{{\mathfrak C}}
\nc{\fD}{{\mathfrak D}}
\nc{\fE}{{\mathfrak E}}
\nc{\fF}{{\mathfrak F}}
\nc{\fG}{{\mathfrak G}}
\nc{\fH}{{\mathfrak H}}
\nc{\fI}{{\mathfrak I}}
\nc{\fJ}{{\mathfrak J}}
\nc{\fK}{{\mathfrak K}}
\nc{\fL}{{\mathfrak L}}
\nc{\fM}{{\mathfrak M}}
\nc{\fN}{{\mathfrak N}}
\nc{\fO}{{\mathfrak O}}
\nc{\fP}{{\mathfrak P}}
\nc{\fQ}{{\mathfrak Q}}
\nc{\fR}{{\mathfrak R}}
\nc{\fS}{{\mathfrak S}}
\nc{\fT}{{\mathfrak T}}
\nc{\fU}{{\mathfrak U}}
\nc{\fV}{{\mathfrak V}}
\nc{\fW}{{\mathfrak W}}
\nc{\fZ}{{\mathfrak Z}}
\nc{\fX}{{\mathfrak X}}
\nc{\fY}{{\mathfrak Y}}
\nc{\fa}{{\mathfrak a}}
\nc{\fb}{{\mathfrak b}}
\nc{\fc}{{\mathfrak c}}
\nc{\fd}{{\mathfrak d}}
\nc{\fe}{{\mathfrak e}}
\nc{\ff}{{\mathfrak f}}
\nc{\fg}{{\mathfrak g}}
\nc{\fh}{{\mathfrak h}}
\nc{\fiI}{{\mathfrak i}}  
\nc{\ffi}{{\mathfrak i}}  
\nc{\fj}{{\mathfrak j}}
\nc{\fk}{{\mathfrak k}}
\nc{\fl}{{\mathfrak{l}}}
\nc{\fm}{{\mathfrak m}}
\nc{\fn}{{\mathfrak n}}
\nc{\fo}{{\mathfrak o}}
\nc{\fp}{{\mathfrak p}}
\nc{\fq}{{\mathfrak q}}
\nc{\fr}{{\mathfrak r}}
\nc{\fs}{{\mathfrak s}}
\nc{\ft}{{\mathfrak t}}
\nc{\fu}{{\mathfrak u}}
\nc{\fv}{{\mathfrak v}}
\nc{\fw}{{\mathfrak w}}
\nc{\fz}{{\mathfrak z}}
\nc{\fx}{{\mathfrak x}}
\nc{\fy}{{\mathfrak y}}

\nc{\al}{{\alpha }}
\nc{\be}{{\beta }}
\nc{\ga}{{\gamma }}
\nc{\de}{{\delta }}
\nc{\del}{{\partial }}
\nc{\vep}{{\varepsilon }}
\nc{\ep}{{\epsilon }}

\nc{\ze}{{\zeta }}
\nc{\et}{{\eta }}
\rc{\th}{{\theta }}

\nc{\vth}{{\vartheta }}

\nc{\io}{{\iota }}
\nc{\ka}{{\kappa }}
\nc{\la}{{\lambda }}
\nc{\vrho}{{\varrho}}
\nc{\si}{{\sigma }}
\nc{\ups}{{\upsilon }}
\nc{\vphi}{{\varphi }}
\nc{\om}{{\omega }}

\nc{\Ga}{{\Gamma }}
\nc{\De}{{\Delta }}
\nc{\nab}{{\nabla}}
\nc{\Th}{{\Theta }}
\nc{\La}{{\Lambda }}
\nc{\Si}{{\Sigma }}
\nc{\Ups}{{\Upsilon }}
\nc{\Om}{{\Omega }}

\nc{\Spec}{\on{Spec}}
\nc{\id}{\on{id}}
\nc{\inv}{ ^{-1}}
\nc{\su}{\subset}
\nc{\ot}{\otimes}
\nc{\un}{\underline}
\nc{\ov}{\overline}
\nc{\uu}{\mathds{1}}

\nc{\shom}{\cH om}
\nc{\Gr}{\cG\tn{r}}

\nc{\se}{\section}
\nc{\sse}{\subsection}
\nc{\ssse}{\subsubsection}

\nc{\md}{\!\!\!\mod}
\nc{\wt}{\widetilde}
\nc{\slc}{\fs \fl_2\bC}
\rc{\sl}{\fs\fl}
\nc{\gr}{\tn{gr}}
\nc{\lan}{\langle}
\nc{\ran}{\rangle}
\nc{\ord}{\tn{ord}}
\nc{\wh}{\widehat}
\nc{\arr}{\ar@<-0.5ex>}
\nc{\re}{\color{red}}
\nc{\bl}{\color{blue}}
\nc{\sst}{\scriptstyle}
\nc{\sss}{\scriptscriptstyle}
\nc{\qchoose}[2]{\left[\!\! \begin{array}{c} #1 \\ #2 \end{array} \!\!\right]}
\nc{\spr}{\tn{spread}}
\nc{\ecd}{\tn{degree}}
\nc{\lce}{\bigg \lceil}
\nc{\rce}{\bigg \rceil}
\nc{\lf}{\lfloor}
\nc{\rf}{\rfloor}

\title{Rank-unimodality of Young's lattice via explicit chain decomposition}

\author{Vivek Dhand}

\begin{document}

\maketitle
\thispagestyle{empty}


\begin{abstract}
Young's lattice $L(m,n)$ consists of partitions having $m$ parts of size at most $n$, ordered by inclusion of the corresponding Ferrers diagrams.  K. O'Hara gave the first constructive proof of the unimodality of the Gaussian polynomials by expressing the underlying ranked set of $L(m,n)$ as a disjoint union of products of centered rank-unimodal subsets.  We construct a finer decomposition which is compatible with the partial order on Young's lattice, at the cost of replacing the cartesian product with a more general poset extension.  As a corollary, we obtain an explicit chain decomposition which exhibits the rank-unimodality of $L(m,n)$.  Moreover, this set of chains is closed under the natural rank-flipping involution given by taking complements of Ferrers diagrams.



\end{abstract}

\se{Introduction}

Young's lattice $L(m,n)$ consists of partitions $\la = (0 \leq \la_1 \leq \dots \leq \la_m \leq n)$, equipped with the following partial order:
\[ \la \leq \la' \iff \la_i \leq \la'_i  \tn{ for all } 1 \leq i \leq m. \]
The {\em rank} of $\la \in L(m,n)$ is defined by $\tn{rk}(\la) = \la_1 + \dots + \la_m$.   The rank generating function of $L(m,n)$ is equal to the Gaussian binomial coefficient:
\[ G(m,n) = \qchoose{m+n}{m}_q = \prod_{i = 1}^m \frac{1-q^{n+i}}{1-q^i} . \]
The unimodality of the coefficients of $G(m,n)$ has been known for a long time, and has been proved using many different techniques (e.g.~see \cite{O,Pr,Z1} for historical background and references).
K. O'Hara gave the first purely {\em combinatorial} proof by expressing the underlying ranked set of $L(m,n)$ as a disjoint union of products of centered rank-unimodal subsets.  More precisely, she defined two statistics on $L(m,n)$ called spread and degree, and decomposed of $L(m,n)$ into centered ranked subsets $U(m,n,s,d)$ consisting of partitions of spread $s$ and degree $d$.  Then she established a rank-preserving bijection:
\[ U(m,n,s,d) \simeq \bigsqcup_{s' < s} U(m-sd,n-2d,s',d') \times U(sn + 2s - 2m,d). \]
Since the product of symmetric unimodal ranked sets is symmetric unimodal, it follows by induction that $U(m,n,s,d)$ and $L(m,n)$ are rank-unimodal.

In this paper, we construct a refinement of spread and degree which we call the signature.  The signature of a partition in $L(m,n)$ is a sequence of non-negative integers $(d_0, \dots, d_k)$ such that $k = \lf n/2 \rf$ and:
\[ m = \sum_{j = 0}^{k} (j+1)d_j . \]
The signature is related to spread and degree by the following formulas:
\[ \spr = d_0 + \dots + d_k \quad \tn{and} \quad \ecd = 1+\min\{0 \leq j \leq k \mid d_j > 0 \}. \]
Let $Q_n(d_0, \dots, d_k)$ denote the subset of $L(m,n)$ consisting of partitions of signature $(d_0, \dots, d_k)$.  We describe a raising and lowering algorithm which provides a covering of $Q_n(d_0, \dots, d_k)$ by saturated chains of length equal to:
\[ \ell_n(d_0, \dots, d_k) = \sum_{j = 0}^k (n - 2j) d_j. \]
This set of chains is stable under the natural rank-flipping involution $\tau$ given by complementation of partitions.  We prove a version of O'Hara's structure theorem which is compatible with Young's partial order.

\sse*{Theorem}  There is a split extension of ranked posets:
\[ \xymatrix{ L(\ell, r) \ar[r] & Q_n(d_0, \dots, d_k)  \ar@<0.4ex>[r] & Q_{n-2r}(d_r, \dots, d_k) \ar@<0.4ex>[l] }  \]
where $\ell = \ell_n(d_0, \dots, d_k)$ and $r = 1 + \min\{0 \leq j \leq k \mid d_j > 0\}$.

By induction, we obtain an explicit $\tau$-stable chain decomposition of $L(m,n)$ so that the chains of a given length can be organized into centered rank-unimodal subposets.

Let us briefly outline the contents of the paper.  In section \ref{signature}, we recall O'Hara's definitions of spread and degree, and we extend them to obtain the signature.  In section \ref{algorithm}, we describe the raising and lowering algorithm and we study the properties of the resulting ``transversal" chains.  We also show that spread, degree, and signature are invariant under this algorithm.  In section \ref{structure}, we prove the above structure theorem for the level sets of the signature map and we construct the explicit chain decomposition which exhibits the rank-unimodality of Young's lattice.

\se{Spread, degree, and signature} \label{signature}

Consider the following poset:
\[ A_n(m) = \{ (a_0, \dots, a_n) \in \bZ_{\geq 0}^{n+1} \mid a_0 + \dots + a_n = m \} \]
where the covering relations are of the form:
\[ (a_0, \dots, a_n) \to (a_0, \dots, a_i -1, a_{i+1} + 1, \dots, a_n). \]
This poset is ranked by the following function:
\[ \tn{rk}(a_0, \dots, a_n) = \sum_{i = 0}^n i a_i. \]
Moreover, there is an isomorphism of ranked posets:
\[ \vphi: L(m,n) \to A_n(m) \]
\[ \vphi(\la) = (a_0, \dots, a_n) \]
where $a_i$ is equal to the number of times $i$ appears in $(\la_1, \dots, \la_m)$.  

Let $\ga: L(m,n) \to L(n,m)$ denote the conjugation isomorphism:
\[ \ga(\la) = \la', \quad  \la' = (0 \leq \la_1' \leq \dots \leq \la_n' \leq m) \]
where $\la_j'$ is equal to the number of $\la_i$ that are greater than or equal to $j$.  Since the Ferrers diagrams of $\la$ and $\la'$ are related by a flip, it follows immediately that $\ga = \ga\inv$.  By composing $\vphi$ and $\ga$, we get an isomorphism:
\[ \psi: L(n,m) \to A_n(m) \]
\[ \psi(\la_1, \dots, \la_n) = (m-\la_n, \la_n - \la_{n-1}, \dots, \la_2 - \la_1, \la_1). \]
In other words, there is a commutative diagram of isomorphisms of ranked posets:
\[  \xymatrix{ L(m,n) \ar@<0.3ex>[rr]^{\ga} \ar[rd]_{\vphi} & & L(n,m) \ar@<0.3ex>[ll]^{\ga}  \ar[ld]^{\psi} \\ & A_n(m) & }  \]
Furthermore, if we define the involution $\tau$ on $A_n(m)$ as follows:
\[ \tau(a_0, \dots, a_n) = (a_n, \dots, a_0), \]
then the isomorphisms in the above diagram also commute with $\tau$.

Recall that O'Hara defines the {\em spread} of $\la \in L(n,m)$ as follows:
\[ \spr(0 \leq \la_1 \leq \dots \leq \la_n \leq m) = \max_{1 \leq i \leq n}(\la_{i + 1} - \la_{i-1}) \]
where $\la_0 = 0$ and $\la_{n+1} = m$.   Let $M(\la)$ denote the set of indices where this maximum value occurs:
\[ M(\la) = \{ 1 \leq i \leq n\mid \la_{i+1} - \la_{i-1} = \spr(\la) \}. \]
Then we have a decomposition:
\[ M(\la) = \bigsqcup_{i} D_i \]
where each $D_i$ is a maximal interval of consecutive integers in $M(\la)$.  O'Hara defines the {\em degree} of $\la \in L(n,m)$ to be:
\[ \ecd(\la) = \sum_{i} \lce \frac{|D_i|}{2} \rce. \]

We transfer these statistics to $A_n(m)$ via the isomorphism $\psi$, i.e. for any $\fa\in A_n(m)$, we define:
\[ \spr(\fa) = \spr(\psi\inv\fa), \quad M(\fa) := M(\psi\inv \fa), \quad \ecd(\fa) = \ecd(\psi \inv \fa). \]
It is straightforward to derive formulas for spread and degree in terms of $A_n(m)$.

\sse{Proposition} If $\fa = (a_0, \dots, a_n) \in A_n(m)$, then:
\[ \spr(\fa) = \max_{0 \leq i \leq n-1}(a_i+ a_{i+1}). \]
\begin{proof}
We know that $\psi\inv\fa = (\la_1, \dots, \la_n)$, where:
\[ \la_i = \sum_{j = n - i + 1}^n a_j. \]
Now, for each $1 \leq i \leq n$, we have:
\[ \la_{i+1} - \la_{i-1} =  \sum_{j = n - i}^n a_j - \sum_{j = n - i + 2}^n a_j  = a_{n-i} + a_{n-i+1}. \]
Therefore:
\[ \spr(\fa) = \max_{1 \leq i \leq n}(\la_{i+1} - \la_{i-1}) = \max_{1 \leq i \leq n}(a_{n-i} + a_{n-i+1}) = \max_{0 \leq i \leq n-1}(a_i + a_{i+1}), \]
where at the last step we have reindexed by $i \mapsto (n-i)$.
\end{proof}  

\sse{Proposition} If $\fa \in A_n(m)$, then:
\[ M(\fa) = \{0 \leq i \leq n-1 \mid a_{i} + a_{i+1} = \spr(\fa) \}. \]

\begin{proof}
If $\la = \psi\inv\fa$, then:
\[ i \in M(\la) \iff \la_{i+1} - \la_{i-1} = \spr(\la) \iff a_{n-i} + a_{n-i+1} = \spr(\fa). \]
Therefore:
\[ M(\fa) = \{1 \leq i \leq n \mid a_{n-i} + a_{n-i+1} = \spr(\fa) \}. \]
Reindexing by $i \mapsto n-i$, we get:
\[ M(\fa) = \{0 \leq i \leq n-1 \mid a_{i} + a_{i+1} = \spr(\fa) \}. \]
\end{proof}

\sse{Remark}  Let $\Ga_n$ denote the path graph with vertex set $\{0, \dots, n\}$ where any two consecutive integers are adjacent.   If $S \su \{0, \dots, n-1\}$, then an {\em edge covering} of $S$ is a collection of pairwise disjoint edges in $\Ga_n$ whose union contains $S$.  The {\em edge covering number} of $S$ is defined to be the number of edges in a minimal edge covering of $S$.   Given $\fa \in A_n(m)$, note that the maximal intervals of consecutive integers in $M(\fa)$ are nothing but the components of $M(\fa)$, thought of as induced subgraph of $\Ga_n$.  If $D$ is any interval of consecutive integers, then $ \lceil |D|/2 \rceil$ is equal to the edge covering number of $D$. Since the components of $M(\fa) $ must have at least one space between them, it follows that $\ecd(\fa)$ is equal to the edge covering number of $M(\fa)$.

A {\em maximal pair} of $\fa = (a_0, \dots, a_n) \in A_n(m)$ is a pair of consecutive entries $(a_i, a_{i+1})$, where $0 \leq i \leq n-1$, such that:
\[ a_i + a_{i+1} = \spr(\fa).  \]
Note that $M(\fa)$ is equal to the set of left indices of maximal pairs of $\fa$.  

We define a natural refinement of spread and degree which we call {\em signature}.  The essential idea is to keep track of the spread and degree as we remove maximal pairs from $\fa \in A_n(m)$.  

Given $\fa \in A_n(m)$, we define the set of {\em active indices} of $\fa$ as follows:
\[ \tn{Act}(\fa) = M(\fa) \cup (1 + M(\fa)). \]
In other words, an index $0 \leq i \leq n$ is active if $(a_i,a_{i+1})$ or $(a_{i-1}, a_i)$ is a maximal pair of $\fa$.  Consider the decomposition of $M(\fa)$ into components:
\[ M(\fa) = \bigsqcup_{i=1}^N [ c_i, c_i + d_i ] \]
where each $ d_i \geq 0$ and $c_{i+1} > c_i + d_i + 1$.  Then:
\[ \tn{Act}(\fa) = \bigsqcup_{i = 1}^N [c_i, c_i  + d_i + 1]. \]
Note that the interval $[c_i, c_i  + d_i + 1]$ has $d_i + 2$ elements. 

If $d_i$ is even, then the part of $\fa$ with indices in $[c_i, c_i  + d_i + 1]$ looks like:
\[ (x,y,x,y, \dots, x,y)  \]
where $x = a_{c_i}$ and $y = a_{c_i + 1}$.

If $d_i$ is odd, then the part of $\fa$ with indices in $[c_i, c_i  + d_i + 1]$ looks like:
\[ (x,y,x,y, \dots, x,y,x)  \]
where $x = a_{c_i}$ and $y = a_{c_i + 1}$.

Let $\om(\fa)$ denote the result of removing the largest possible number of maximal pairs from $\fa$.  Following the above discussion, we see that it does not matter in what order we remove the maximal pairs;  in the end $\om(\fa)$ will contain all entries $a_j$ such that $j \notin \tn{Act}(\fa)$, along with those entries $a_{c_i}$ where $d_i$ is odd.  The number of maximal pairs removed while calculating $\om(\fa)$ is equal to $\ecd(\fa)$, and $\spr(\fa)$ is equal to the sum of the entries of each maximal pair, so:
\[ \om(\fa) \in A_{n - 2r}(m - rs) \]
where $r = \ecd(\fa)$ and $s = \spr(\fa)$. Also note that $\spr(\om(\fa)) < \spr(\fa)$ by construction.

\sse{Proposition} The maps $\om$ and $\tau$ commute.
\begin{proof}
Let $\fa \in A_n(m)$.  Note that $i \in M(\fa)$ if and only if $n-i-1 \in M(\tau \fa)$.  Similarly, $i \in \tn{Act}(\fa)$ if and only if $n-i \in \tn{Act}(\tau\fa)$.  It follows immediately that $\om(\tau \fa)= \tau \om(\fa)$.
\end{proof}

Let $k = \lfloor n/2\rfloor$.  We define the {\em signature}:
\[ \si : A_n(m) \to \bZ_{\geq 0}^{k + 1} \]
by induction on $n$.  If $n \leq 1$,  let $\si(\fa) = m$.  If $n > 1$, let:
\[ \si(\fa) = (0^{\ecd(\fa) - 1} , \spr(\fa) - \spr(\om(\fa)), \si(\om(\fa) )). \]

\sse{Proposition} Let $\fa \in A_n(m)$ with $m > 0$.  If  $\si(\fa) = (d_0, \dots, d_k)$, then:
\[ \ecd(\fa) = 1 + \min \{ {0 \leq j \leq k} \mid d_{j} > 0 \}, \quad \spr(\fa) = d_0 + \dots + d_k . \]

\begin{proof}
Let $r = \ecd(\fa)$ and $s = \spr(\fa)$.  We know that:
\[ d_{r-1} = s - \spr(\om(\fa)) > 0 \quad \tn{and} \quad d_j = 0 \tn{ for } 0 \leq j \leq r-2, \]
so:
\[ 1 + \min \{ {0 \leq j \leq k} \mid d_{j} > 0 \} = 1 + (r - 1) = r \]
which proves the first equation.  Since $\si(\om(\fa)) = (d_r, \dots, d_k)$, we have:
\[  \spr(\om(\fa)) = d_r + \dots + d_k \]
by induction.  Therefore:
\[ s = d_{r-1} + d_r + \dots + d_k = d_0 + \dots + d_k \]
which proves the second equation.
\end{proof}

We are now ready to use the signature to refine O'Hara's decomposition.  Given $n, d_0, \dots, d_k \geq 0$, we define:
\[ Q_n(d_0, \dots, d_k) = \{ \fa \in A_n(m) \mid   \si(\fa) = (d_0, \dots, d_k) \} \]
where:
\[ m = \sum_{j = 0}^k (j+1) d_j. \]

\sse{Proposition} $Q_n(d_0, \dots, d_k)$ is stable under $\tau$.

\begin{proof}
Let $\fa = (a_0, \dots, a_n) \in Q_n(d_0, \dots, d_k)$.  Then $\tau \fa = (a_n, \dots, a_0)$, and:
\[ \spr(\tau \fa) = \spr(\fa) \quad \tn{and} \quad \ecd(\tau \fa) = \ecd(\fa). \]
Therefore:
\begin{eqnarray*}
\si(\tau\fa) &=& (0^{\ecd(\tau\fa) - 1} , \spr(\tau\fa) - \spr(\om(\tau\fa)), \si(\om(\tau\fa) ))\\
  &=& (0^{\ecd(\fa) - 1} , \spr(\fa) - \spr(\tau\om(\fa)), \si(\tau\om(\fa) )) \\
  &=& (0^{\ecd(\fa) - 1} , \spr(\fa) - \spr(\om(\fa)), \si(\om(\fa) )) \\
  &=& \si(\fa) 
\end{eqnarray*}
where we have used that $\om$ and $\tau$ commute and that $\si(\tau\om(\fa)) = \si(\om(\fa))$ by induction.
\end{proof}

In particular, $Q_n(d_0, \dots, d_k)$ is a centered subposet of $A_n(m)$.

\sse{Lemma}  Let $n\geq 0$, $k = \lfloor n/2\rfloor$, and  $d_0, \dots, d_k \geq 0$.  

(1) The restriction of $\om$ defines a surjective map:
\[  \om_r: Q_n(d_0, \dots, d_k) \to Q_{n - 2r}(d_r, \dots, d_k ) \]
where $r = 1 + \min\{0 \leq j \leq k \mid d_j > 0\}$.

(2) There is a $\tau$-stable decomposition:
\[ Q_n(d_0, \dots, d_k) = \bigsqcup_{\fb \in Q_{n-2r}(d_r, \dots, d_k)} \om_r\inv \fb. \]

\begin{proof}
(1) Since $\si(\fa) = (0^{r-1}, \spr(\fa) - \spr(\om(\fa), \si(\om(\fa))) = (d_0, \dots, d_k)$, we see that:
\[ d_0 = \dots = d_{r-2} = 0, \quad d_{r-1} > 0, \quad \tn{and} \quad \si(\om(\fa)) = (d_r, \dots, d_k). \]
Therefore, $\om$ restricts to a map:
\[  \om_r : Q_n(d_0, \dots, d_k) \to Q_{n - 2r}(d_r, \dots, d_k ). \]
To prove that $\om_r$ is surjective, we define a right inverse:
\[  \be_r : Q_{n-2r}(d_r, \dots, d_k) \to Q_{n}(d_0, \dots, d_k ). \]
\[ \fb \mapsto ((s,0)^r, \fb) \]
where $s = d_0 + \dots + d_k$.  Since $d_{r-1} > 0$, we have:
\[ \spr(\fb) = d_r + \dots + d_k < s. \]
Therefore, $M(\be_r(\fb)) = \{0, \dots, 2r-1\}$, which means that:
\[ \ecd(\be_r(\fb)) = r \quad \tn{and} \quad  \si(\be_r(\fb)) = (0^{r-1}, d_{r-1}, d_r , \dots, d_r). \]
It follows that $\tn{Act}(\be_r(\fb)) = \{0, \dots, 2r\}$ and $\om_r(\be_r(\fb)) = \fb$.


(2) This follows immediately from (1) and the fact that $\om_r$ commutes with $\tau$.
\end{proof}

\sse{Remark} Signature is the natural ``completion" of spread and degree.  More precisely, to any $\fa \in A_n(m)$, Conca \cite{Co} associates a certain tableau consisting of $m$ boxes with entries in $\{0, \dots, n\}$.  The spread is equal to the size of the longest column and the degree is equal to the size of the shortest row.  Therefore, spread and degree together describe a rectangular block inside the tableau.  We  can remove this block and calculate the spread and degree of what remains, and so on.  In the end, we will get a complete description of the shape of the tableau, corresponding to the signature.  Indeed, if $\si(\fa) = (d_0, \dots, d_k)$, then the associated tableau has $d_j$ rows of size $(j+1)$.   We can also define the signature in terms of certain tropical polynomials derived from the secant varieties of the rational normal curve in projective space \cite{Dh}.  In this way, we see that there is a natural geometric interpretation of spread, degree, and signature. 

\se{The raising and lowering algorithm} \label{algorithm}

We define a weight function on $A_n(m)$ as follows:
\[ \tn{wt}: A_n(m) \to \bZ \]
\[ \tn{wt}(\fa) = mn - 2 \tn{rk}(\fa) \]
Note that weights of all elements of $A_n(m)$ are congruent to $mn$ mod 2.  Also note that, for any $\fa \in A_n(m)$:
\[ \tn{wt}(\tau \fa) = mn - 2 \sum_{i = 0}^n (n-i) a_{i} = 2\sum_{i = 0}^n i a_i - mn = - \tn{wt}(\fa) \]
which implies that the weights of $A_n(m)$ are distributed symmetrically around zero.  From now on, when working with elements of $A_n(m)$ we will refer to their weights instead of their ranks.  Note that the elements of highest weight correspond to the elements of lowest rank and vice versa.


We say that $\fa\in A_n(m)$ is {\em initial} (resp.  {\em terminal}) if $(a_0,0)$ (resp. $(0,a_n)$) is a maximal pair.
Note that it is possible for an element to be both initial and terminal, e.g. $\fa = (2,0,1,0,0,2)$.

We are now ready to describe the raising and lowering algorithm.

\sse*{The raising algorithm}  Let $(a_i,a_{i+1})$ be a maximal pair of $\fa \in A_n(m)$.  If $i \geq 1$, then $a_{i+1} \geq a_{i-1}$ because $a_i + a_{i+1} \geq a_{i-1} + a_i$.  If $a_{i+1} > a_{i-1}$, then we decrement $a_{i+1}$ and increment $a_i$.   If $a_{i+1} = a_{i-1}$, then we start over with the maximal pair $(a_{i-1}, a_i)$.  If the current maximal pair is $(a_0, a_1)$, then we decrement $a_1$ and increment $a_0$.  Once $a_1 = 0$, we end the chain.  Note that this endpoint is an initial element of $A_n(m)$.

\sse*{The lowering algorithm}  Let $(a_{i-1},a_i)$ be a maximal pair of $\fa \in A_n(m)$.   If $i \leq n-1$, then $a_{i-1} \geq a_{i+1}$ because $a_{i-1} + a_i \geq a_i + a_{i+1}$.  If $a_{i-1} > a_{i+1}$, then we decrement $a_{i-1}$ and increment $a_i$.   If $a_{i-1} = a_{i+1}$, then we start over with the maximal pair $(a_i, a_{i+1})$.  If the current maximal pair is $(a_{n-1}, a_n)$, then we decrement $a_{n-1}$ and increment $a_n$.  Once $a_{n-1} = 0$, we end the chain.  Note that this endpoint is a terminal element of $A_n(m)$.

\sse{Proposition}  Spread and degree are invariant under the raising and lowering algorithm.

\begin{proof}
We will prove the statement for the raising algorithm. It is straightforward to adapt the proof to the lowering algorithm.  Let $\fa = (a_0, \dots, a_n) \in A_n(m)$ and let $(a_i,a_{i+1})$ be a maximal pair of $\fa$.  We may assume that $a_{i-1} < a_{i+1}$.  In particular, $(a_{i-1},a_i)$ is not a maximal pair of $\fa$.  Applying the raising algorithm, we obtain the element:
\[ \fa' = (a'_0, \dots, a'_n) =  (a_0, \dots, a_{i} + 1, a_{i+1} - 1, \dots, a_n). \]
Now let us calculate the spread of $\fa'$.  There are exactly two sums of consecutive entries which are affected:
\[  a'_{i-1} + a'_{i} = a_{i-1} + a_i + 1 \quad\tn{and}\quad  a'_{i+1} + a'_{i+2} = a_{i+1} + a_{i+2} - 1. \]
Even if $(a_{i+1},a_{i+2})$ was a maximal pair, we see that $(a'_{i+1},a'_{i+2})$ is not a maximal pair and does not matter for calculating the spread.  On the other hand, if $(a'_{i-1},a'_i)$ is a maximal pair, then:
\[  a'_{i-1} + a'_{i} = a_{i-1} + a_i + 1 \leq a_i + a_{i+1}. \]
It follows that $(a'_i, a'_{i+1}) = (a_i + 1, a_{i+1} - 1)$ is a maximal pair of $\fa'$, and $\spr(\fa') = \spr(\fa)$.   

Now let us compare $M(\fa)$ to $M(\fa')$.  Note that $i$ lies in both $M(\fa)$ and $M(\fa')$ but $i-1 \notin M(\fa)$ and $i+1 \notin M(\fa')$.  Therefore, the maximal consecutive interval in $M(\fa)$ which contains $i$ is of the form $[i, i+d]$ for some $d \geq 0$ and the maximal consecutive interval in $M(\fa')$ which contains $i$ is of the form $[i-d', i]$ for some $d' \geq 0$.  The change from $M(\fa)$ to $M(\fa')$ is:
\[ [i-d', i-2] \sqcup [i,i+d] \mapsto [i-d', i] \sqcup [i+2,i+d] \]
where $[i+2,i+d]$ (resp.~ $[i-d',i-2]$) is empty if $d < 2$ (resp.~ $d' < 2$).
In the degree calculation, the sum of the relevant terms remains constant:
\[ \lce \frac{d'-1}{2} \rce +  \lce \frac{d+1}{2} \rce =  \lce \frac{d'+1}{2} \rce +  \lce \frac{d-1}{2} \rce.  \]
Therefore, $\ecd(\fa) = \ecd(\fa')$.
\end{proof}

\sse{Proposition} The signature is invariant under the raising and lowering algorithm.

\begin{proof}
Let $\fa = (a_0, \dots, a_n) \in A_n(m)$ and let $(a_i,a_{i+1})$ be a maximal pair of $\fa$ such that $a_{i-1} < a_{i+1}$.  Applying the raising algorithm, we obtain the element:
\[ \fa' = (a'_0, \dots, a'_n) =  (a_0, \dots, a_{i} + 1, a_{i+1} - 1, \dots, a_n). \]
Since spread and degree are invariant under the raising and lowering algorithm, it suffices to show that $\om(\fa) = \om(\fa')$.  

Note that $i \in M(\fa) \cap M(\fa')$.  In particular, we have $i, i+1 \in \tn{Act}(\fa) \cap \tn{Act}(\fa')$.  Also, $i$ is the leftmost element of its component in $M(\fa)$ and the rightmost element of its component in $M(\fa')$.  In either case, we may assume that $i$ and $i+1$ are removed while calculating $\om(\fa)$ and $\om(\fa')$.  Therefore, if $\tn{Act}(\fa) = \tn{Act}(\fa')$, then $\om(\fa) = \om(\fa')$.

There are two possible cases where $\tn{Act}(\fa) \neq \tn{Act}(\fa')$.

(1) If $i+2 \notin M(\fa)$ and $i+1 \in M(\fa)$, then $i+2 \in \tn{Act}(\fa)$ but $i+2 \notin \tn{Act}(\fa')$.  

In this case, $i+2$ survives in $\om(\fa')$ because it is not active in $\fa'$.  Since $i$ is the leftmost endpoint of its component in $M(\fa)$, and $i+1$ is the rightmost, we see that $i+2$ survives in $\om (\fa)$ as well.  Therefore $\om(\fa) = \om(\fa')$.

(2) If $i-2 \notin M(\fa)$ and $i-1 \in M(\fa')$, then $i-1 \notin \tn{Act}(\fa)$ but $i-1 \in \tn{Act}(\fa')$.

In this case, $i-1$ survives in $\om(\fa)$ because it is not active in $\fa$.  Note that $i$ is the rightmost endpoint of its component in $M(\fa')$, and $i-1$ is the leftmost, and $i+1$ is active in $\fa'$.  Therefore, either $i-1$ or $i+1$ survives in $\om(\fa')$.  By assumption:
\[ a'_{i-1} = a_{i-1} = a_{i+1} -1 = a'_{i+1}, \]
so we conclude that $\om(\fa) = \om(\fa')$.
\end{proof}

Given $\fa$ in $A_n(m)$ and $i \in M(\fa)$, we can apply the raising and lowering algorithm starting at the maximal pair $(a_i,a_{i+1})$ and obtain a {\em transversal chain} $T_i(\fa)$ whose highest (resp.~lowest) weight element is initial (resp.~terminal).

\sse{Proposition} Let $\fa \in A_n(m)$.  There is a bijection between set of transversal chains containing $\fa$ and the set of components of $M(\fa)$.

\begin{proof}
Let $i \in M(\fa)$.  Suppose $i - 1 \in M(\fa)$ and start the raising algorithm at the maximal pair $(a_i, a_{i+1})$.  Since $a_{i+1} = a_{i-1}$, the algorithm will simply restart with the maximal pair $(a_{i-1}, a_i)$.  In other words, the raising algorithm affects the entries of $\fa$ if and only if $i$ is the leftmost endpoint of its component in $M(\fa)$.  Similarly, the lowering algorithm affects the entries of $\fa$ if and only if $i$ is the rightmost endpoint of its component in $M(\fa)$.  It follows that, if $i$ and $j$ lie in the same component of $M(\fa)$, then $T_i(\fa) = T_j(\fa)$.

It remains to show that if $i$ and $j$ are in distinct components of $M(\fa)$, then $T_i(\fa) \neq T_j(\fa)$.  Indeed, let $i_1$ (resp.~$j_1$) denote the leftmost element of the component containing $i$ (resp.~$j$).   Then the color of the edge of $T_i(\fa)$ (resp.~$T_j(\fa)$) which enters $\fa$ from above is $i_1$ (resp.~$j_1$).  Since $i_1 \neq j_1$, the transversal chains are distinct.
\end{proof}

Recall that the Hasse diagram of $A_n(m)$ has an edge-coloring where the $i$-th color corresponds to the covering relation:
\[ (a_0, \dots, a_n) \to (a_0, \dots, a_{i-1} -1, a_{i} + 1, \dots, a_n). \]
Any saturated chain in $A_n(m)$ is uniquely determined by its highest weight element along with the sequence of colors obtained by following the chain from highest to lowest weight in the Hasse diagram.   One of the key properties of transversal chains is that their color sequences are non-decreasing.

\sse{Proposition}
If $\fa = (a_0,0,a_2, \dots, a_n) \in A_n(m)$ is an initial element, then the color sequence of $T_0(\fa)$ is:
\[ (1^{a_0 - a_2},2^{a_0 - a_2 - a_3}, \dots, j^{a_0 - a_j - a_{j+1}} ,\dots, n^{a_0 - a_n} ). \]
Furthermore, the terminal element of $T_0(\fa)$ is:
\[ \fb = (a_2, \dots, a_n, 0, a_0) = \tau (a_0, 0, \tau(a_2, \dots, a_n)) \quad \tn{and} \quad T_0(\fa) = T_{n-1}(\fb). \]
\begin{proof}
We calculate $T_0(\fa)$ by applying the lowering algorithm to the leftmost maximal pair in $\fa$, namely $(a_0,0)$.  Since $a_0 \geq a_2$, we decrement $a_0$ and increment $a_1$ exactly $(a_0-a_2)$ times:
\[ (a_0, 0, a_2, \dots, a_n) \mapsto \dots \mapsto (a_2, a_0 - a_2, a_2, \dots, a_n) \]
which means the color sequence starts with $1^{a_0 - a_2}$.  Then we consider the maximal pair $(a_0 - a_2, a_2)$ and compare $a_0 - a_2$ to $a_3$.  We must have $a_0 - a_2 \geq a_3$ since $a_0 \geq a_2 + a_3$.  Therefore, we decrement $(a_0 - a_2)$ and increment $a_2$ exactly $(a_0 - a_2 - a_3)$ times:
\[ (a_2, a_0 - a_2, a_2, \dots, a_n) \mapsto \dots \mapsto (a_2, a_3, a_0  - a_3, a_3, \dots, a_n),  \]
hence appending $2^{a_0-a_2-a_3}$ to the color sequence.  Continuing in this way, we see that the $j$-th color will appear $(a_0 - a_j - a_{j+1})$ times in the color sequence, where $a_{n+1} = 0$ by definition.  

To find the termial element of $T_0(\fa)$, note that we are essentially moving $a_0$ to the right while shifting $a_i$ two spots to the left for each $i \geq 2$.  Therefore, the terminal element in this chain will be $\fb = (a_2, \dots, a_n, 0, a_0)$.  

Now we calculate $T_{n-1}(\fb)$ by applying the raising algorithm to the rightmost maximal pair in $\fb$, namely $(0,a_0)$.  Since $a_0 \geq a_n$, we decrement $a_0$ and increment $0$ exactly $(a_0 - a_n)$ times:
\[ (a_2, \dots, a_n, 0, a_0) \mapsto \dots \mapsto (a_2, \dots, a_n, a_0 - a_n, a_n) \]
which means the reverse of the color sequence starts with $n^{a_0 - a_n}$.  Then we consider the maximal pair $(a_n, a_0 - a_n)$ and compare $a_0 - a_n$ to $a_{n-1}$.  We must have $a_0 - a_n \geq a_{n-1}$ since $a_0 \geq a_n + a_{n-1}$.  Therefore, we decrement $(a_0 - a_n)$ and increment $a_n$ exactly $(a_0 - a_{n} - a_{n-1})$ times:
\[ (a_2, \dots, a_n, a_0 - a_n, a_n) \mapsto \dots \mapsto (a_2, \dots, a_{n-1}, a_{0} - a_{n-1}, a_{n-1}, a_n), \]
hence appending $(n-1)^{a_0 - a_n - a_{n-1}}$ to the reverse color sequence.  Continuing in this way, we see that the $(n+1-j)$-th color will appear $(a_0 - a_{n+2-j} - a_{n+1-j})$ times, where $a_1 = 0$.  In this way, we obtain the revese of the color sequence of $T_0(\fa)$, so $T_0(\fa) = T_{n-1}(\fb)$.
\end{proof}


Next we show that the set of transversal chains is stable under $\tau$.  

\sse{Proposition}  If  $\fa \in A_n(m)$ is an initial element and $\fb$ is the terminal element of $T_0(\fa)$, then:
\[ \tau T_0(\fa) = T_0(\tau \fb).  \]

\begin{proof}
Note that $\tau$ swaps initial and terminal elements and replaces the $j$-th color with $(n + 1 -j)$-th color.  Therefore, the initial element of $\tau T_0(\fa)$ is $\tau \fb$ and its color sequence is the reverse of the color sequence of $T_0(\fa)$ with the $j$-th color replaced by $(n + 1 - j)$:
\[ (1^{a_0 - a_n}, 2^{a_0 - a_{n-1} - a_{n}}, \dots,(n + 1 - j)^{a_0 - a_j - a_{j+1}},\dots, (n-1)^{a_0 - a_2 - a_3}, n^{a_0 - a_2}). \]
On the other hand, since $\fb = \tau (a_0, 0, \tau(a_2, \dots, a_n))$, we see that 
\[ \tau \fb = (a_0, 0, a_n, a_{n-1}, \dots, a_2)\] 
and therefore the color sequence of $T_0(\tau\fb)$ is:
\[ ( 1^{a_0 - a_n} , 2^{a_0 - a_{n-1} - a_{n}}, \dots, j^{a_0 - a_{n + 1 - j}- a_{n+2-j}} ,\dots, (n-1)^{a_0 - a_2 - a_3}, n^{a_0 - a_2}). \]
Since the two chains have the same initial element and color sequence, we conclude that $ \tau T_0(\fa) = T_0(\tau \fb)$.
\end{proof}

We finish this section by deriving two more useful properties of $Q_n(d_0, \dots, d_k)$.

\sse{Proposition} $Q_n(d_0, \dots, d_k)$ has a unique element of highest weight, namely $\fh = (h_0, \dots, h_n)$, where $h_i = 0$ if $i$ is odd and:
\[ h_{2i} = \sum_{j = i}^k d_j \]
for $0 \leq i \leq k$.

\begin{proof}
First let's check that $\fh \in Q_n(d_0, \dots, d_k)$.  Note that $h_{2i} - h_{2i+2} = d_i$.  By construction, $\fh$ has the property that $h_0 \geq h_2 \geq \dots \geq h_{2k}$, which implies that:
\[ h_{2i} = \spr(h_{2i}, \dots, h_n) \]
for each $0 \leq i \leq k$.  In particular:
\[ \om(\fh) = (h_{2r}, \dots, h_n)\]
where:
\[ r = 1 + \min\{ 0 \leq i \leq k \mid h_{2i} > h_{2i+2} \} = 1 + \min\{0 \leq i \leq k \mid d_i > 0\}.\]  
Therefore:
\[ \si(\fh) = (0^{r-1}, h_{2r-2} - h_{2r}, \si(\om(\fh))) = (0^{r-1}, d_{r-1}, d_r, \dots, d_k) \]
where we have used that $\si(\om(\fh)) = (d_r, \dots, d_k)$ by induction.

Now suppose that $\fm = (m_0, \dots, m_n)$ is a highest weight element of $Q_n(d_0, \dots, d_k)$.  In particular, $\fm$ is initial, so $m_1 = 0$ and $m_0 = \spr(\fm) = d_0 + \dots + d_k$.  

Now consider the element $\fm' = (m_2, \dots, m_n) \in Q_{n-2}(d_1, \dots, d_k)$.  If $\fm'$ is not an initial element of $Q_{n-2}(d_1, \dots, d_k)$, then the raising algorithm will produce an element $\fn' \in Q_n(d_1, \dots, d_k)$ such that $\tn{wt}(\fn') = \tn{wt}(\fm') + 2$.  

But then the element $\fn = (m_0, 0, \fn') \in Q_n(d_0, \dots, d_k)$ and $\tn{wt}(\fn) = \tn{wt}(\fm) + 2$, which contradicts our assumption.  Therefore, $\fm'$ is initial in $Q_n(d_1, \dots, d_k)$.  It follows that $m_3 = 0$ and $m_2 = \spr(\fm') = d_1 + \dots + d_k$.  By repeating this argument for each $(m_{2i}, \dots, m_n)$ for $1 \leq i \leq k$, we see that $\fm = \fh$. 
\end{proof}

Since $Q_n(d_0, \dots, d_k)$ is defined as a level set for the signature map, it follows immediately that it has a covering by transversal chains, and it turns out that they all have the same length.

\sse{Proposition} All the transversal chains in $Q_n(d_0, \dots, d_k)$ have length equal to:
\[ \ell_n(d_0, \dots, d_k) = \sum_{j = r-1}^k (n-2j)d_j \]
where $r = 1 + \min\{0 \leq j \leq k \mid d_j > 0\}$.
\begin{proof}
Let $\fa \in Q_n(d_0, \dots, d_k)$ be an initial element and let $T_0(\fa)$ be the corresponding transversal chain.  By definition, $a_1= 0$ and $0 \in M(\fa)$, so:
\[  a_0 = \spr(\fa) = d_0 + \dots + d_k = d_{r-1} + \dots + d_k = s.\]
Then the color sequence of $T_0(\fa)$ is:
\[ (1^{s - a_2},2^{s - a_2 - a_3}, \dots, j^{s - a_j - a_{j+1}} ,\dots, n^{s - a_n} ), \]
and therefore its length is:
\[ ns - 2(a_2 + \dots + a_n). \]
Now:
\[ a_0 + \dots + a_n = d_0 + 2d_1 + \dots + (k+1)d_k  = m\]
so:
\[ a_2 + \dots + a_n = m - a_0 =  m - s, \]
so the length of $T_0(\fa)$ depends only on $n,d_0, \dots, d_k$.  Let $\ell_n(d_0, \dots, d_k)$ denote the length of each transversal chain in $Q_n(d_0, \dots, d_k)$. 

To calculate $\ell_n(d_0, \dots, d_k)$, it suffices to calculate the length of $T_0(\fh)$ for the unique highest weight element $\fh = (h_0, \dots, h_n) \in Q_n(d_0, \dots, d_k)$, where $h_i = 0$ for $i$ odd and 
\[ h_{2i} = \sum_{j = i}^k d_j \]
for $0 \leq i \leq k$.  The length of $T_0(\fh)$ is equal to:
\begin{eqnarray*}
nh_0 - 2\sum_{i = 1}^k h_{2i} &=& n h_0 - 2 \sum_{i = 1}^k \sum_{j = i}^k d_j  \\
&=&n h_0 - 2 \sum_{j = 1}^k \sum_{i = 1}^j d_j  \\
&=&n \sum_{j = 0}^k d_j - 2 \sum_{j = 1}^k j d_j  \\
&=&  \sum_{j = 1}^k (n-2j) d_j.
\end{eqnarray*}
The desired formula now follows from the fact that $d_j = 0$ for $0\leq j \leq r-2 $.
\end{proof}

\se{The structure theorem} \label{structure}

Let $F, E, $ and $B$ be ranked posets.  We say that $E$ is a {\em split extension} of $B$ by $F$ if there exists a diagram:
\[ \xymatrix{ F \ar[r] & E \ar@<0.4ex>[r]^\om & B \ar@<0.4ex>[l]^\be } \]
such that:

(1) $\om$ is a surjective order-preserving map.

(2) $\be$ is a order-preserving section of $\om$.

(3) For each $b \in B$, there is an isomorphism of ranked posets $\al_b: F \simeq \om\inv b$. 

It follows that there exists a covering of $E$ of the form: 
\[ E = \bigsqcup_{b \in B} F_b, \]
where $\al_b: F \simeq F_b$ and, if $b_1$ covers $b_2$ in $B$, then: 
\[ \tn{rk}(\al_{b_1}(f)) = \tn{rk}(\al_{b_2}(f)) + 1\] 
for any $f \in F$. 

We are now ready to prove our main structure theorem for Young's lattice order.

\sse{Theorem}  There exists a split extension of ranked posets:
\[ \xymatrix{ L(r,\ell) \ar[r] & Q_{n}(d_0, \dots, d_k) \ar@<0.4ex>[r]^<<<<<{\om_r} & Q_{n-2r}(d_r, \dots, d_k) \ar@<0.4ex>[l]^<<<<<{\be_r} } \]
where $\ell = \ell_n(d_0, \dots, d_k)$ and $r = 1 + \min\{0 \leq j \leq k \mid d_j > 0\}$.
\begin{proof}
We have already proved that there is a decomposition:
\[ Q_n(d_0, \dots, d_k) = \bigsqcup_{\fb \in Q_{n-2r}(d_r, \dots, d_k)} \om_r\inv(\fb) \]
with an order-preserving section $\be_r: Q_{n-2r}(d_r, \dots, d_k) \to Q_{n}(d_0, \dots, d_k)$.  Therefore, it remains to show that $\om_r$ is order-preserving and that there is an isomorphism of ranked posets:
\[ \om_r\inv(\fb) \simeq L(r, \ell) \]
for each $\fb\in Q_{n-2r}(d_r, \dots, d_k)$.

Let us prove the latter statement first.  We recursively define a map:
\[ \De: \om_r\inv(\fb) \to L(r,\ell) \]
\[ \fa  \mapsto (\de_1, \dots, \de_r). \]
Let $s = d_0 + \dots + d_k$.  We define $\de_1$ to be the number of steps it takes for the raising algorithm to reach an initial element $(s,0,\fa')$ starting with the leftmost maximal pair of $\fa$. Then we define $(\de_2, \dots, \de_r) = \De(\fa')$.  

We need to check that $0 \leq \de_1 \leq \dots \leq \de_r \leq \ell$.  Since $\ell$ is equal to the length of each transversal chain in $Q_{n-2j}(d_{j}, \dots, d_k)$, for each $0 \leq j \leq r-1$, we see that each $\de_i \leq \ell$.   

By definition, the first entry of $\De(s,0,\fa')$ is zero.  By induction, $\de_2  \leq \dots \leq \de_r$, so it suffices to show that $\de_1 \leq \de_2$.  
It turns out there a nice formula for $\de_1$. Indeed, if $(a_i,a_{i+1})$ is the leftmost maximal pair of $\fa$, then:
\[ \de_1 = (i+1)s - a_{i} - 2(a_0 + \dots + a_{i-1}) =  a_{i+1} + is - 2(a_0 + \dots + a_{i-1}).\]
Also, 
\[ \fa' = (a_0, \dots, a_{i-1}, a_{i+2}, \dots ) \]
If $(a_j, a_{j+1})$ is the leftmost maximal pair of $\fa'$, then $j \geq i+2$ and:
\[ \de_2 = a_{j+1} + js - 2(a_0 + \dots + a_{i-1} + a_{i+2} + \dots + a_{j-1}). \]
Therefore:
\begin{eqnarray*}
 \de_2 - \de_1 &=& a_{j+1} + a_i + (j-i-1)s - 2(a_{i+2} + \dots + a_{j-1}) \\
 &\geq& a_{j+1} + a_i + (a_{i+1} + 2a_{i+2} + \dots + 2a_{j-1} + a_{j}) - 2(a_{i+2} + \dots + a_{j-1}) \\
 &=& a_{j+1} + a_i + a_{i+1} + a_{j} = 2s \geq 0.
 \end{eqnarray*}





Now we need to define an inverse for $\De$ and prove that both $\De$ and $\De\inv$ are order-preserving for Young's lattice order.

Given $\la = (\la_1, \dots, \la_r) \in L(r,\ell)$, let $\fb_1(\la)$ be the element obtained by following $T_0(\be_r(\fb))$ for $\la_r$ steps.  For $2 \leq i \leq r$, define $\fb_i(\la)$ to be the element obtained by following $T_0(\fb_{i-1}(\la))$ for $\la_{r+1-i}$ steps.   Note that each $\fb_i(\la)$ must lie in $\om_r\inv(\fb)$ because each $\la_i \leq \ell$ and $\om$ is invariant under the raising and lowering algorithm.  Now we define:
\[ \De\inv: L(r,\ell) \to \om_r\inv(\fb) \]
\[ (\la_1, \dots, \la_r) \mapsto \fb_r(\la). \]
By construction, $\De$ and $\De\inv$ are inverse to each other.

To show that $\De$ is order-preserving, let $\fp, \fq \in \om_r\inv(\fb)$ such that $\fp$ covers $\fq$.  We claim that $\De(\fp)$ covers $\De(\fq)$ in $L(r,\ell)$.  Indeed, $\fp$ is obtained from $\fq$ by one step of the lowering algorithm at a maximal pair $(q_i, q_{i+1})$.

If  $(q_i, q_{i+1})$ is the leftmost maximal pair of $\fq$, then it is also the leftmost maximal pair of $\fp$.  Therefore, $\fq$ and $\fp$ lie on the same transversal chain, separated by one step, and so $\de_1(\fp) = \de_1(\fq) + 1$.   Since the raising algorithm will send both $\fp$ and $\fq$ to the same initial element, it follows that $\De(\fp)$ and $\De(\fq)$ agree in all subsequent entries.  Therefore, $\De(\fp)$ covers $\De(\fq)$ in $L(r,\ell)$.

On the other hand, if $(q_i, q_{i+1})$ is not the leftmost maximal pair of $\fq$, then the raising algorithm follows the same steps for $\fp$ and $\fq$, which yields initial elements $(s,0,\fp')$ and $(s,0,\fq')$ such that $\fp'$ covers $\fq'$.  By induction, $\De(\fp')$ covers $\De(\fq')$ in $L(r-1,\ell)$.

To show that $\De\inv$ is order-preserving, let $\la, \la' \in L(r,\ell)$ such that $\la$ covers $\la'$.  Then $\la_i = \la_i' + 1$ for some $1 \leq i \leq r$ and $\la_j = \la_j'$ for all $j \neq i$.  




If $i = 1$, then $\fb_p(\la) = \fb_p(\la')$ for $1 \leq p \leq r-1$.  It follows that $\fb_r(\la)$ and $\fb_r(\la')$ are separated by one step on the transversal chain: 
\[ T_0(\fb_{r-1}(\la)) = T_0(\fb_{r-1}(\la')). \] 
Therefore, $\De\inv(\la)$ covers $\De\inv(\la')$ in $\om_r\inv(\fb)$.  

On the other hand, if $i > 1$, then $\fb_p(\la) = \fb_p(\la')$ for $1 \leq p \leq r-i$. It follows that $\fb_{r-i+1}(\la)$ and $\fb_{r-i+1}(\la')$ lie on the same transversal chain, separated by one step.  By assumption, $(\la_1, \dots, \la_{i-1}) = (\la'_1, \dots, \la'_{i-1})$.  We claim that $\fb_{r-i+2}(\la)$ covers $\fb_{r-i+2}(\la')$.  Indeed, we know that:
\[  \fb_{r-i+1}(\la) = (s,0,\fp) \quad\tn{and}\quad  \fb_{r-i+1}(\la') = (s,0,\fq) \]
where $\fp$ is obtained from $\fq$ by applying one step of the lowering algorithm.  Furthermore, $\fp$ (resp.~$\fq$) is obtained by applying the lowering algorithm $\la_i$ times (resp.~$\la'_i$ times) to $\fb_{r-i}(\la) = \fb_{r-i}(\la')$.  Since $\la_{i-1} = \la_{i-1}' \leq \la'_i = \la_i - 1$, we see that $T_0(s,0,\fp)$ and $T_0(s,0,\fq)$ have the same color sequence for at least $\la_{i-1}$ steps, and therefore $\fb_{r-i+2}(\la)$ covers $\fb_{r-i+2}(\la')$.  By induction, we conclude that $\De\inv(\la)$ covers $\De\inv(\la')$ in $\om_r\inv(\fb)$. 

Finally, we show that $\om_r$ is order-preserving.  In fact, we will prove a stronger result.  Suppose $\fp$ covers $\fq$ in $Q_{n}(d_0, \dots, d_k)$ and that they do not lie on the same transversal chain.  Let $(s,0,\fp')$ (resp.~$(s,0,\fq')$) denote the initial element obtained by  applying the raising algorithm to the leftmost maximal pair of $\fp$ (resp.~$\fq$).  We will prove that $\fp'$ covers $\fq'$ in $Q_{n-2}(d_1,\dots, d_k)$. 

Let $(q_j,q_{j+1})$ denote the leftmost maximal pair of $\fq$.  If we start the raising algorithm at this maximal pair, then we end up with the inital element:
\[ (s, 0, q_0, \dots, q_{j-1}, q_{j+2},\dots, q_n). \]
Now suppose $\fp$ and $\fq$ are related by an edge of color $(i+1)$, so they only differ in two places: $p_{i} = q_{i} - 1$ and $p_{i+1} = q_{i+1} + 1$. 

If $0 \leq i \leq j-2$, then:
\[ \fp' = (q_0, \dots , q_i - 1, q_{i+1} + 1, \dots, q_{j-1}, q_{j+2}, \dots, q_n) \]
so $\fp'$ covers $\fq'$.

If $i = j-1$, then $p_{j} +  p_{j+1} = (q_j +1) + q_{j+1} > s$, which is a contradiction.

If $i = j$, then $\fp' = \fq'$, so $\fp$ and $\fq$ lie on the same transversal chain, which contradicts our assumption.

If $i = j + 1$, then $(p_{j}, p_{j+1}) = (q_j, q_{j+1} - 1)$ is not a maximal pair of $\fp$.  If $(q_{j+1}, q_{j+2})$ is not a maximal pair of $\fq$, then $\ecd(\fp) < \ecd(\fq)$, which is a contradiction.   On the other hand, if $(q_{j+1}, q_{j+2})$ is a maximal pair of $\fq$, then $(p_{j+1},p_{j+2})$ now becomes the leftmost maximal pair of $\fp$.  Then one step of the raising algorithm applied to $(p_{j+1},p_{j+2}) = (q_{j+1} - 1, q_{j+2}+1)$ gives us $\fq$ back again, which means that $\fp$ and $\fq$ lie on the same transversal chain, contradicting our assumption. 

Finally, if $j + 2 \leq i \leq n-1$, then $(q_j,q_{j+1})$ is the leftmost maximal pair of both $\fp$ and $\fq$, and so:
\[ \fp' = (q_0, \dots, q_{j-1},q_{j+2}, \dots, q_{i}-1, q_{i+1} +1, \dots, q_n) \]
so $\fp'$ covers $\fq'$.
\end{proof}

\sse{Remark} It follows immediately from the structure theorem that $Q_n(d_0, \dots, d_k)$ and $L(m,n)$ are rank-unimodal.  Indeed, if $d_0 > 0$, then $r = 1$, and:
\[ Q_n(d_0, \dots, d_k)  = \bigsqcup_{\fb \in Q_{n-2}(d_1, \dots, d_k)} T_0(\be_1(\fb)). \]
Therefore, $Q_{n}(d_0, \dots, d_k)$ is covered by the disjoint union of saturated chains of the same length whose highest weight elements form a rank-unimodal set, which means that $Q_{n}(d_0, \dots, d_k)$ is rank-unimodal.  On the other hand, if $d_0 = 0$, then $r \geq 2$, and $Q_n(d_0, \dots, d_k)$ is a split extension of a rank-unimodal poset by $L(\ell,r)$.  Now $L(\ell,r)\simeq A_r(\ell)$ has a decomposition into centered rank-unimodal subposets of the form $Q_r(d'_0, \dots, d'_t)$, where $t = \lfloor r/2 \rfloor$.  
It follows that there is a decomposition of $Q_n(d_0, \dots, d_k)$ into centered rank-unimodal subposets isomorphic to a split extension of $Q_{n-2r}(d_0, \dots, d_k)$ by $Q_r(d'_0, \dots, d'_t)$.  By induction, each of these split extensions has a $\tau$-stable chain decomposition where the highest weight elements of chains of a given length form a rank-unimodal set.

\sse{Remark} Note that our proof of the rank-unimodality of Young's lattice does not use the fact that $\om_r$ is order-preserving.  Indeed, we have only used the weaker property that the highest weight elements of the chains of a given length form a rank unimodal set.  However, we have included this stronger property to make the statement of the structure theorem more elegant, and it seems to be an important fact in its own right.

 \end{document}